# Exact Solution of Abel Differential Equation with Arbitrary Nonlinear Coefficients

Ali Bakhshandeh Rostami

Federal University of Rio de Janeiro (UFRJ), Rio de Janeiro, Brazil

**Abstract**

This paper is dedicated to present an exact solution for a nonlinear differential equation so-called Abel equation. This equation was known as one of the group of unsolvable differential equations. The present method is applicable for any arbitrary form of nonlinear coefficients of Abel equation. On the other hands, the presented solution in this paper is verified with exact solution for some restricted forms of Abel equation that has been reported by Polyanin and Zaytsev [5]. Also, numerical solution is utilized to verify this method. All verifications have approved that the present method is strongly reliable and exact to solve Abel equation, analytically.

**Keywords**: Abel Equation of second Kind, Exact Solution, Nonlinear Equation

**Highlights:**

- Exact solution of a nonlinear equation so-called Abel Equation is presented in this paper.
- Abel Equation was known as an unsolvable equation.

## 1  Introduction

Two main categories as numerical and analytical methods exist to analyze the mathematical problems. Numerical methods are suitable for general problems, but the computational procedure is complex, while in most cases, the efficiency and accuracy of the results are relatively lower than those of the analytical methods [1]. Therefore, among the existing methods, the most efficient is analytical method. Thus, many scholars have applied various analytical methods to find a reliable solution for different mathematical problems [2]. On the other hands, for many decades, nonlinear differential equations such as Able equation were and are in center of attentions to find an applicable and reliable analytical solution for them. Abel equation appears to reduce the order of many higher orders nonlinear problems. Hence, every day, the exact solution is demanded for Abel differential equation in first and second kind.

The Abel differential equation frequently found in the modeling of real problems in varied areas such as big picture modeling in oceanic circulation [8], in problems of magneto-statics [9] and [10], control theory [11], cosmology [12], Fluid Mechanics [13] and [14], Solid Mechanics[19], biology [17] and cancer therapy [18].

So far, the great number of works has attempted to propose a solution for both kind of Abel equation. For instance, Bougofa [6] reduced the general equation of second kind to a canonical form and proposed a particular construction that can solve this equation. Moreover, Mak and Harko [15] have reported a method to generate the general solution of Abel equation of the first kind from a particular solution. The Gülsu et al [16] have utilized an approximate method by using the shifted Chebyshev expansion. Therefore, Abel differential equation is approximately transformed to a system of nonlinear equations that gives a solvable system in matrix form. Panayotounkos [4] has used a certain type of parameterization to achieve an analytical solution. The author has found that the solution is not reliable in his presented form and require to careful scrutiny. It is worth to mention that in all these references above, the restricted form of analytical solutions are given. Since, after an accurate revise on the last paper, finally the present paper is reporting a reliable exact solution for Abel differential equation while as we know, the general solution to Abel differential equations was an open problem so far.

## 2  Abel Equation

The Abel equation of second kind is any nonlinear ordinary differential equation that is quadratic in unknown functions and is written in general form as:

$$[g_1(x)\, U + g_0(x)]U'_x = f_2(x)\, U^2 + f_1(x)\, U + f_0(x). \qquad (1)$$

Here the notations $U'_x = dU/dx$, $U''_x = d^2U/dx^2$, … are used for the total derivatives. Also, g1, g0, f2, f1 and f0 are known as variable and arbitrary form of coefficients of the Abel equation. The first kind of Abel equation has cubic fraction in unknown function as follow:

$$U'_x = f_3(x) U^3 + f_2(x) U^2 + f_1(x) U + f_0(x). \tag{2}$$

The equation (1) and also restricted form of (2) can be reduced to normal form (or canonical form) by using various admissible functional transformations as mentioned in [3] and [5]. The normal form of Abel equation is expressed as,

$$yy'_x - y = Q(x). \tag{3}$$

In which y and Q are sufficiently differentiable functions in the subinterval $[x_1, x_2]$ and Q(x) is called non-homogeneous part. Polyanin and Zaytsev [5] have reported some exact analytical solutions of equation (3) for special forms of inhomogeneity ($Q(x)$). In present study, the author attempts to construct an exact analytic solution of (3) for general form of $Q(x)$ in the above mentioned subinterval by using the Bougoffa construction (Appendix A). This means that the general form of Abel equation (in the second kind (1) and also restricted form of the first kind (2)) admits an exact analytic solution.

## 3 Solution of Normal form

For solving the normal form of the Abel equation (3), the functional transformation is introduced as [4]

$$y(x) = h(x) V(\xi(x)). \tag{4}$$

By inserting Equation (4) into Equation (3), the canonical form of the Abel equation is rewritten as:

$$h^2 \xi'_x V V'_\xi + h^2 h'_x V^2 - hV = Q(x). \tag{5}$$

Where h(x), ξ(x) and V(ξ(x)) need to be determined. Now, by introducing a new differentiable function as k(x), the last equation is rewritten as:

$$(h^2 \xi'_x V + k) V'_\xi - 2Q(x) = (-h^2 \xi'_x V + k) V'_\xi - 2h^2 h'_x V^2 + 2hV. \tag{6}$$

This equation can be split into two Abel equations as follow. The variable coefficients of Abel equation have been indicated in following equations.

$$\left( \underbrace{h^2 \xi'_x V}_{g_{1a}} + \underbrace{k}_{g_{0a}} \right) V'_\xi = \underbrace{G(x) + 2Q(x)}_{f_{0a}} \qquad (a),$$

$$\left( \underbrace{-h^2 \xi'_x V}_{g_{1b}} + \underbrace{k}_{g_{0b}} \right) V'_\xi = \underbrace{2h^2 h'_x}_{f_{2b}} V^2 - \underbrace{2h}_{f_{1b}} V + \underbrace{G(x)}_{f_{0b}} \qquad (b). \tag{7}$$

Here, G(x) is a differentiable function that shall be concluded. By taking into account $\xi(x) = x$ for functional transformation of Equation (4) and after proper integrations, the Bougoffa construction coefficients (Equation (61)) for two Abel equations of (7) (a and b) is calculated as:

$$B_{1a}(x) = 1, \; B_{2a}(x) = 1, \; \lambda_a = \lambda = \frac{2\,k(x)}{h^2(x)\xi'_x} \quad (a),$$

$$B_{1b}(x) = h^4(x), B_{2b}(x) = e^{\frac{4}{\lambda}\int \frac{dx}{h(x)}}, \; \lambda_b = -\frac{2\,h^4(x)k(x)}{h^2(x)e^{\frac{4}{\lambda}\int\frac{dx}{h(x)}}} \quad (b).$$

(8)

Now, Bougoffa construction is formed,

$$V^2 + \lambda V - 2\int \frac{G(x) + 2Q(x)}{h^2(x)}\,dx = 0 \quad (a),$$

$$V^2 - \lambda V + \frac{2}{h^4(x)}\int h^2(x)\,G(x)\,dx = 0 \quad (b).$$

(9)

From now on and for brevity, we use G instead G(x), Q instead of Q(x) and so on. In set of equations of (9), the subsidiary functions G and h, as well as V have to be determined. Therefore, the determination of all subsidiary functions is discussed in the following sections.

### 3.1  Determination of h(x)

The subsidiary function h can be found by the third items of Equation (8)-b as follow:

$$h(x) = \frac{1}{\lambda}(x + \varphi). \tag{10}$$

That φ is integration constant. Also, from third item of Equation (8)-b can be deduced that $\lambda_b = -\lambda^{-3}$ whereas $\lambda_a = \lambda$ that λ is unknown constant. The author uses of the Julia construction (Appendix B) to determine λ. By applying the Julia construction on Abel equation of (7)-a and then, put it equals to the Bougoffa construction for the same Abel equation,

$$h(x) = \frac{1}{2}(x + \varphi), \quad k(x) = h(x)^2, \quad \lambda = 2. \tag{11}$$

### 3.2  Determination of V(x)

For determining the V(x), the author uses the real roots of the quadratic equations (a and b) of (9) as follow:

$$V = \frac{-\lambda \pm \sqrt{\lambda^2 + 8\int \frac{G + 2Q}{h^2}}}{2} \quad (a), \tag{12}$$

$$V = \frac{\lambda \pm \sqrt{\lambda^2 - \frac{8}{h^4} \int h^2 G}}{2} \qquad (b).$$

By equating a and b of (12), one equation obtains as follow:

$$\sqrt{\lambda^2 h^4 - R} = h^2 \sqrt{\lambda^2 + K} - 2h^2 \lambda. \tag{13}$$

Where,

$$K(x) = 8 \int \frac{G + 2Q}{h^2} dx, \quad \text{and} \quad R(x) = 8 \int h^2 G \, dx. \tag{14}$$

The square of Equation (13) is,

$$3\lambda^2 h^4 - 4h^4 \lambda \sqrt{\lambda^2 + K} + h^4(\lambda^2 + K) + R = 0. \tag{15}$$

Differentiate from (15) and substitute derivative of K and R by expressions of (14), results the following equation:

$$(\lambda^2 + K)^{\frac{3}{2}} - 4\lambda(\lambda^2 + K) + \left(3\lambda^2 + \frac{4(G + Q)}{h'_x h}\right)\sqrt{\lambda^2 + K} - \frac{4\lambda(G + 2Q)}{h'_x h} = 0. \tag{16}$$

The last equation is a cubic one of $\sqrt{\lambda^2 + K}$. Now, by substituting of Equation (11) into (16) and after some manipulations, the cubic equation is rewritten as:

$$(1 + K)^{\frac{3}{2}} - 4(1 + K) + \left(3 + \frac{4(G + Q)}{x + \varphi}\right)\sqrt{1 + K} - \frac{4(G + 2Q)}{x + \varphi} = 0. \tag{17}$$

By the well-known transformation,

$$\sqrt{1 + K} = Z + \frac{4}{3}. \tag{18}$$

Equation (17) reduces to the cardano form as,

$$Z^3 + p Z + q = 0. \tag{19}$$

Where

$$\begin{aligned} p &= -\frac{a^2}{3} + b, & q &= 2\left(\frac{a}{3}\right)^3 - \frac{ab}{3} + c, \\ a &= -4, & b &= \left(3 - c - \frac{2Q}{h}\right), & c &= -\frac{2(G + 2Q)}{h}. \end{aligned} \tag{20}$$

It is well-known that the solution of the cubic equation (19) can be expressed in exact analytic form, depends on the sign of discriminant

$$D = \left(\frac{p}{3}\right)^3 + \left(\frac{q}{2}\right)^2. \tag{21}$$

Thus, the roots of (19) are the following:

Case 1: $D < 0$ and $p < 0$

$$Z_1 = 2\sqrt{-\frac{p}{3}} \, Cos \, \frac{\beta}{3}, \qquad Z_2 = -2\sqrt{-\frac{p}{3}} \, Cos \, \frac{\beta - \pi}{3},$$

$$Z_3 = -2\sqrt{-\frac{p}{3}} \, Cos \, \frac{\beta + \pi}{3}, \quad Cos \, \beta = -\frac{q}{2\sqrt{-\left(\frac{p}{3}\right)^3}} \quad 0 < \beta < \pi. \tag{22}$$

Case 2: $D > 0$

$$Z_1 = \sqrt[3]{-\frac{q}{2} + \sqrt{D}} + \sqrt[3]{-\frac{q}{2} - \sqrt{D}},$$

$$Z_2 = -0.5\left(\sqrt[3]{-\frac{q}{2} + \sqrt{D}} + \sqrt[3]{-\frac{q}{2} - \sqrt{D}}\right) + i\frac{\sqrt{3}}{2}\left(\sqrt[3]{-\frac{q}{2} + \sqrt{D}} - \sqrt[3]{-\frac{q}{2} - \sqrt{D}}\right),$$

$$Z_3 = -0.5\left(\sqrt[3]{-\frac{q}{2} + \sqrt{D}} + \sqrt[3]{-\frac{q}{2} - \sqrt{D}}\right) - i\frac{\sqrt{3}}{2}\left(\sqrt[3]{-\frac{q}{2} + \sqrt{D}} - \sqrt[3]{-\frac{q}{2} - \sqrt{D}}\right). \tag{23}$$

Case 3: $D = 0$

$$Z_1 = 2\sqrt[3]{-\frac{q}{2}}, \quad Z_2 = Z_3 = -\sqrt[3]{-\frac{q}{2}}. \tag{24}$$

By substituting Equation (18) into (12)-a, the function V(x) is concluded as follow:

$$V = Z + \frac{1}{3}. \tag{25}$$

In the last equation, V is a function of unknown parameter Z that is concluded by Equations (22) to (24). On the other side, it can be deduced from Equation (20) that the parameter Z is in terms of the unknown function G(x) by means of coefficient c. Hence, to calculate the V(x), it is necessary to determine the coefficient c of (20). Eventually, the solution of Abel normal form (3) is written in terms of the unknown coefficient c as follows:

$$y(x) = \frac{1}{2}(x + \varphi)\left[Z(x) + \frac{1}{3}\right]. \tag{26}$$

That Z(x) is calculated from Equation (19).

### 3.3 Determination of Coefficient c

The coefficient c is defined as:

$$c = -\frac{2(G + 2Q)}{h}. \tag{27}$$

To determine the above coefficient, the author starts rewriting the Abel equations (7)-a and (7)-b as follow:

$$\begin{aligned}-h^2\xi'_x VV'_\xi &= h^2 V'_\xi - G - 2Q & (a), \\ -h^2\xi'_x VV'_\xi &= 2h^2 h'_x V^2 - 2h V + G - h^2 V'_\xi & (b).\end{aligned} \tag{28}$$

The following Riccati equation has been obtained as a consequence of equating of (28)-a and (28)-b.

$$V'_x = \frac{h'_x}{h} V^2 - \frac{1}{h} V + \frac{G + Q}{h^2}, \tag{29}$$

The Riccati equation can be reduced to normal form by the following functional transformation [5, part 0.1.4-6],

$$\omega(x) = \frac{2hh'_x V + h'_x + h'^2_x}{h'^2_x}. \tag{30}$$

After manipulations, the normal form of Riccati equation (31) is expressed as:

$$\omega'_x = \omega^2 + \frac{h'^2_x}{4h^2}\left[\frac{4(G+Q)}{hh'_x} - \frac{1}{h'^2_x} - \frac{2}{h'_x} + 1\right]. \tag{31}$$

The following set of equations satisfies (31).

$$\omega_p = Z + \frac{1}{3} \quad \text{and} \quad \omega'_{p_x} = \frac{(G + 2Q)}{h\left[Z + \frac{4}{3}\right]}. \tag{32}$$

By [5, part 0.1.4-3], it is known that, if relations of (32) constitute a particular solution of (31), then the general solution is written as:

$$\begin{aligned}\omega_g &= Z + \frac{1}{3} + \tau(x), \quad \text{that} \quad \tau(x) = \frac{\Phi(x)}{C_3 - \int \Phi(x)dx}, & (a), \\ \Phi(x) &= Exp\left(\int 2\left(Z + \frac{1}{3}\right) dx\right) & (b).\end{aligned} \tag{33}$$

where $C_3$ is an integration constant. From the cubic equation of (17), it has been verified that the particular solution of the Riccati equation ($V_p = Z + 1/3$) satisfies the Abel equations (7). Hence, the Abel equation has to be satisfied by the general solution of the Riccati equation, as well. Therefore, it is interpreted as $\tau(x)$ will vanish. This condition is formulated as:

$$\lim_{x \to \pm\infty} \omega_g = \lim_{x \to \pm\infty} \left(Z + \frac{1}{3}\right) \quad or \quad \lim_{x \to \pm\infty} \tau(x) \to 0 \tag{34}$$

The Equation (34) is defined by improper integral as follow:

$$\int_{-\infty}^{+\infty} \tau'_x = 0 \tag{35}$$

Where

$$\tau'_x = \frac{\Phi'_x}{C_3 - \int \Phi(x)dx} + \frac{\Phi^2}{[C_3 - \int \Phi(x)dx]^2} \tag{36}$$

Thus:

$$\int_{-\infty}^{+\infty} \frac{\Phi'_x}{C_3 - \int \Phi(x)dx} + \int_{-\infty}^{+\infty} \frac{\Phi^2}{[C_3 - \int \Phi(x)dx]^2} = 0 \tag{37}$$

The last equation is solved if each integral goes to zero. Hence, the following expressions can be deduced:

$$\int_{-\infty}^{+\infty} \Phi'_x = 0 \tag{1}$$

$$C_3 - \int_{-\infty}^{x} \Phi(x)dx \neq 0 \quad and \quad \neq singular \tag{2} \tag{38}$$

$$\int_{-\infty}^{+\infty} \left|\frac{\Phi}{C_3 - \int \Phi(x)dx}\right| \to 0 \tag{3}$$

where the function $\Phi$ has to be determined. In the following, some mathematical logical statements are presented that help to evaluate this function. From Appendix C and also condition 1 simultaneously, we understand that the derivative of $\Phi$ ($\Phi'_x$) can be an odd function because the integration of an odd function in symmetrical interval is zero.

Also, if suppose $C_3=0$ then from condition 2, it deduces that the function $\Phi$ should be non-zero and non-singular function. Based on Appendix C and condition 2 (if $C_3=0$ and x goes to infinity) simultaneously, the function $\Phi$ should not be an odd function. On the other side, Condition 3 gives one more condition for $\Phi$ that this function is more convenient to be constant function or at least converge on a constant value at infinity.

In brief, the functions $\Phi'_x(x)$ and $\Phi(x)$ are more convenient to suppose as odd and semi-constant even functions, respectively. With regards to former conclusion about $\Phi'_x(x)$ and $\Phi(x)$, the function $\Phi'_x(x)$ can be proposed as follow:

$$\Phi'_x = \frac{\sin x}{|x|} \quad where \quad x \to \pm\infty \tag{39}$$

That is an odd function. From [7, type 3.757],

$$\int_0^{+\infty} \frac{\sin x}{x} = \frac{\pi}{2} \tag{40}$$

Then,

$$\int_{-\infty}^{+\infty} \frac{\sin x}{|x|} = 0 \tag{41}$$

On the other side, the function $\Phi(x)$ is calculated by following relation

$$\Phi(x) = \int \frac{\sin x}{|x|} = sgn(x)\, Si(x) \tag{42}$$

The right hand side of (42) is an semi-constant even function that sgn(x) represents the sign of x. Moreover, the Si(x) gives the Sine integral that is defined as [7, type 8.232]:

$$Si(x) = \int_0^z \frac{sint}{t} dt \quad or \quad Si(x) = -\frac{\pi}{2} + \sum_{k=1}^{\infty} \frac{(-1)^{k+1} x^{2k-1}}{(2k-1)(2k-1)!} \tag{43}$$

The Si(x) is an entire function of x with no branch cut discontinuities and is semi-constant function that converges approximately on 1.5 at infinity.

The defined functions $\Phi'_x(x)$ and $\Phi(x)$, in (39) and (42), satisfy all conditions that discussed in foregoing paragraph. Now, by equating (33)-b and (42) and taking a natural logarithm (Ln) of both sides,

$$\int 2\left(Z + \frac{1}{3}\right) dx = Ln[sgn(x)\, Si(x)] \tag{44}$$

Derivative of (44) becomes:

$$Z + \frac{1}{3} = \frac{Sin\ x}{2\ sgn(x)\ |x|\ Si(x)} \tag{45}$$

Thus, Equation (33)-a is rewritten:

$$\omega_g = \frac{Sin\ x}{2\ sgn(x)\ |x|\ Si(x)} + \frac{sgn(x)\ Si(x)}{C_3 - \int sgn(x)\ Si(x)\ dx} \tag{46}$$

The last expression obeys the assertion (34), $\lim_{x \to \pm\infty} \tau(x) = 0$. By differentiating (45) and equate with second of (32):

$$\frac{G + 2Q}{h(Z + 4/3)} = \frac{Cos\ x\ [|x|Sgn(x)Si(x)] - Sin\ x\ [Si(x) + Sin\ x]}{2\ x^2\ Si^2(x)} \tag{47}$$

And finally, from (27) and for all range of x:

$$c = \frac{0.5\ \psi\ Sin\ 2x - \left[2 + \frac{1}{sgn(x)|x|}\right]\psi\ Sin^2 x + 2\psi^2 \left[Cos\ x - \frac{Sin\ x}{sgn(x)|x|}\right] - Sin^3 x}{-2\ \psi^3} \tag{48}$$

where,

$$\psi = sgn(x)\ |x|\ Si(x) \tag{49}$$

The equations (48) are reliable for entire range of x except x=0. Also, the effect of high order of inhomogeneity is considered in the next section.

## 4  Effect of order of inhomogeneity

Function Q(x) (non-homogenous part of normal form) with $n^{th}$ order is defined as:

$$Q(x) = [f(x)]^n \tag{50}$$

where n is the highest degree of polynomial that has been made up by base function of f(x). The base function f(x) includes any arbitrary function such as exponential, trigonometry, logarithm and so on. For inserting the effect of order of Q(x), equation (26) requires modifying as follow:

$$y(x) = \frac{1}{2}\ (x + \varphi) \left[\beta\ Z(x) + \frac{1}{3}\right] \tag{51}$$

where β is called factor of order effect of inhomogeneity. This factor has been defined in Table **1** for different type of base function.

**Table 1:** Relations for factor of order effect of inhomogeneity as $Q(x) = [f(x)]^n$ for different types of the base function

| Base Function f(x) | β |
|---|---|
| x | -0.252 Ln (n)+0.991 |
| Logarithm e.g. Ln(x) | -0.0053 n²+0.0425n+0.981 |
| Exponential e.g. $e^x$ | -0.12 Ln (n)+0.3218 |
| Trigonometry e.g. Sin(x) | 1 |

These relations of β are also applicable for rational type of inhomogeneity as $Q(x) = \frac{[f(x)]^m}{[f(x)]^k}$. Note that, in rational function,

$$n = m - k \quad \rightarrow \quad \begin{cases} n \geq 2 & \rightarrow \quad \beta = from\ Table\ 1\ based\ on\ f(x) \\ n < 2 & \rightarrow \quad \beta = 1 \end{cases} \quad (52)$$

If Q(x) consists of combination of different types of f(x), relation of β is selected from **Table 1** based on f(x) with greatest degree.

## 5  Verification and discussion

A new method has been presented in section 3 and 4 for solving the Abel equation with any arbitrary coefficients. In the following section, the author is going to verify the presented method with numerical method and reported exact solutions by Polyanin and Zaytsev [5] for special functions. The numerical method that is employed in verifying is "Adams-Bashforth method" that is a Predictor-Corrector method. This method has been considered by Lomax et al [20] comprehensively.

Polyanin and Zaytsev [5] have reported exact solution of canonical form of Abel equation in their section 1.3 for some limited forms of inhomogeneity (Q(x)). The different type of non-homogenous functions and their solutions have been selected in three categories as polynomial, rational and exponential functions.

### 5.1  Polynomial Inhomogeneity

Two type of polynomial functions as linear ($Q(x) = Ax + B$- equation No. 1.3.1.2 of [5]) and cubic ($Q(x) = Sx + 2Ax^2 + 2A^2x^3$- equation No. 1.3.1.14 of [5]) are investigated in this section so that A, B and S are arbitrary parameters. As shown in Figure 1 and Figure 2, the present method solves the normal form of Abel equation with polynomial functions of different orders accurately.

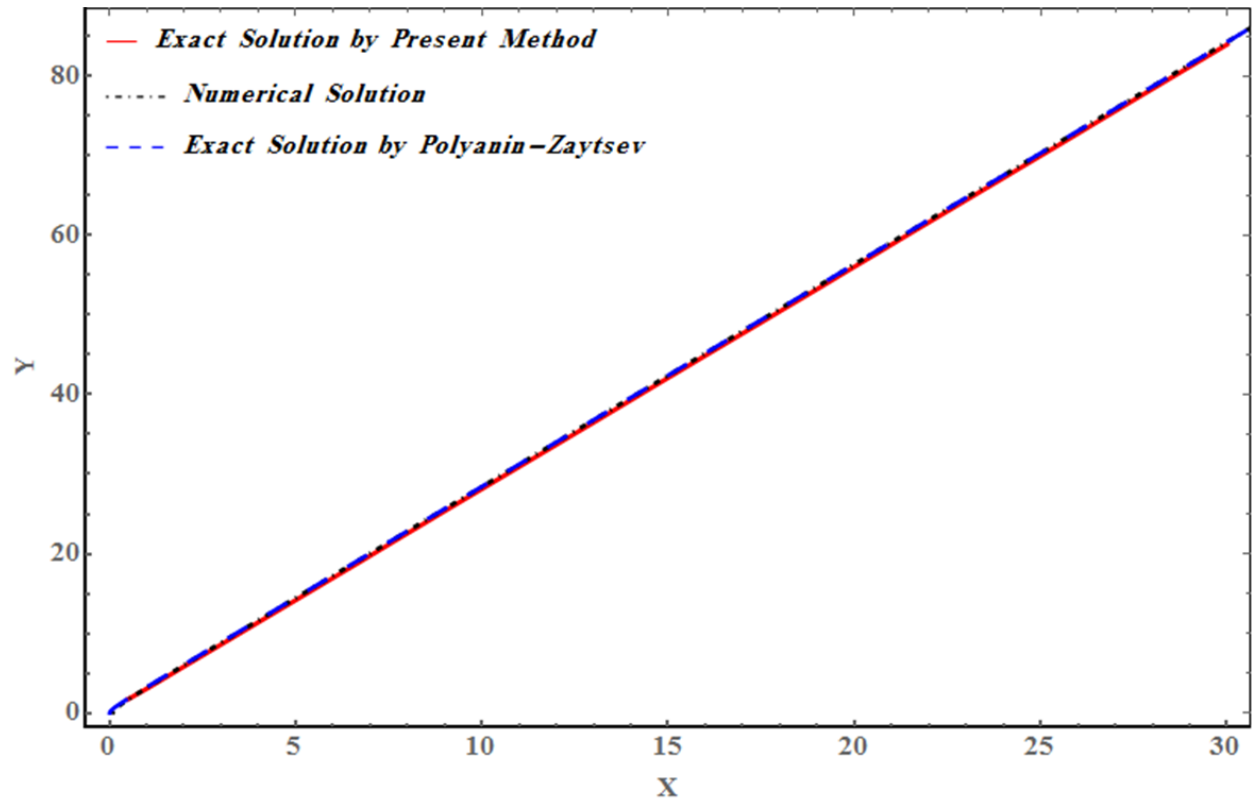

**Figure 1**. Solution of $yy'_x - y = Ax + B$ that A and B are arbitrary parameters. The Solution of present method (with β=1 in equation (51)) has been compared by numerical solution and also exact solution of the same function that reported by Polyanin-Zaytsev [5] (Equation 1.3.1.2).

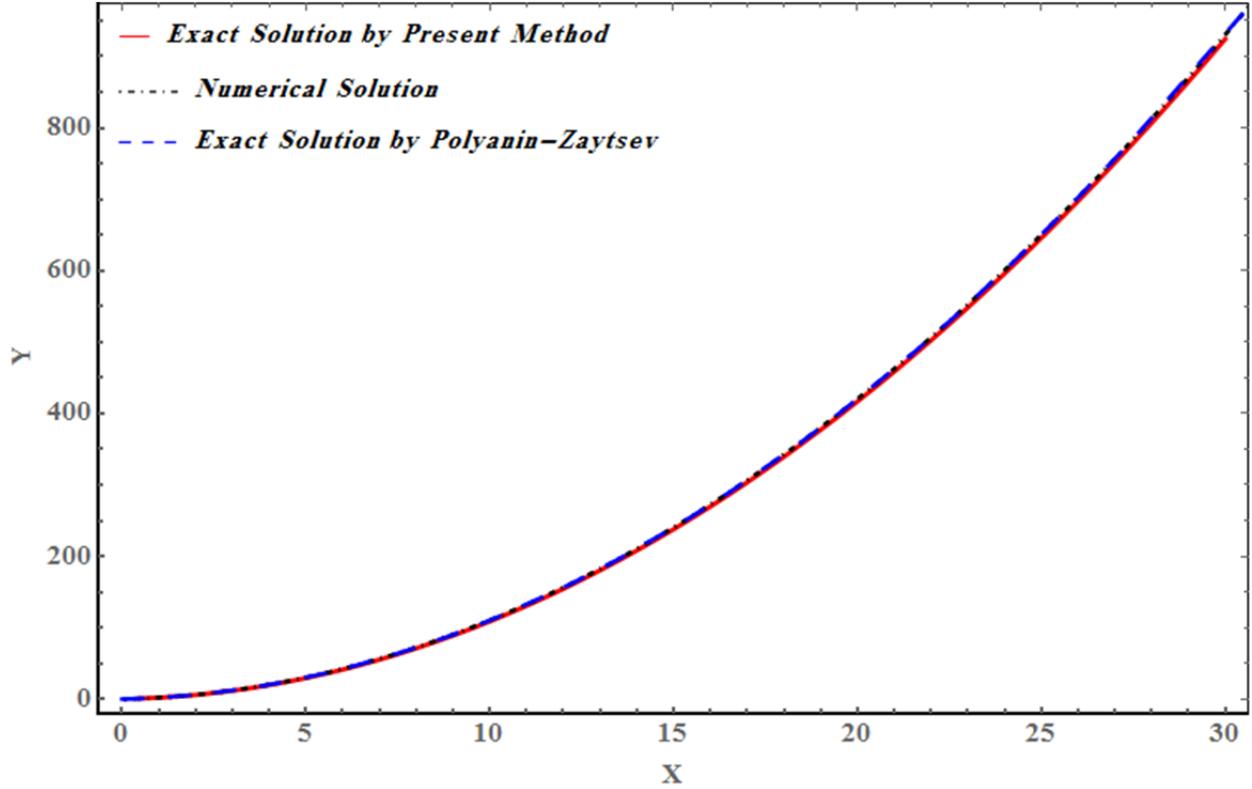

**Figure 2.** The comparison of results are calculated by the present method(β in equation (51) is 0.705), the numerical method and Polyanin-Zatsev's exact solution ([5]-Equation 1.3.1.14). The equation is $yy'_x - y = \frac{4}{9}x + 2Ax^2 + 2A^2x^3$ that A is arbitrary parameter.

## 5.2 Rational Inhomogeneity

According to Polyanin-Zaytsev´s exact solutions for different type of rational functions, two functions are given in this part. The first function is $Q(x) = \frac{A}{x} + \frac{A^2}{x^3}$ that corresponds to Equation No. 1.3.1.7 in [5]. The solution of this function has been presented in parametric form in [5] as follow:

$$x = \sqrt{2A}\tau^{-1}\sqrt{\tau - Ln(1+\tau) - C}$$
$$y = \sqrt{2A}\left[\frac{1+\tau}{\tau}\sqrt{\tau - Ln(1+\tau) - C} - \frac{\tau}{2\sqrt{\tau - Ln(1+\tau) - C}}\right] \quad (53)$$

where C is integration constant. The set of equation (53) has some typos and requires modifying as:

$$x = \sqrt{2A}\tau^{-1}\sqrt{\tau + Ln(1+\tau) - C}$$
$$y = \sqrt{2A}\left[\frac{1+\tau}{\tau}\sqrt{\tau + Ln(1+\tau) - C} - \frac{\tau}{2\sqrt{\tau + Ln(1+\tau) - C}}\right] \quad (54)$$

The modified version of last solution has been presented with results of numerical method and the present method in Figure 3. According to this figure, exact solution by present method shows good accuracy.

The second function that is considered here is equation No. of 1.3.1.33 of [5] as $Q(x) = \frac{A}{x^2}$. The exact solution of this function is reported by polyanin-Zaytsev [5] as:

$$U_1 = \frac{\tau}{2}\left[I_{-\frac{2}{3}}(\tau) + I_{\frac{4}{3}}(\tau)\right] + \frac{1}{3}I_{\frac{1}{3}}(\tau), \quad U_2 = U_1^2 + \tau^2\left[I_{\frac{1}{3}}(\tau)\right]^2, \quad U_3 = \frac{2}{3}\tau^2\left[I_{\frac{1}{3}}(\tau)\right]^2 - 2U_1U_2$$
$$x = -2\sqrt[3]{\frac{A}{36}}\tau^{\frac{4}{3}}U_2^{-1}\left[I_{\frac{1}{3}}(\tau)\right]^2, \quad y = -3\sqrt[3]{\frac{A}{36}}\tau^{-\frac{2}{3}}U_2^{-1}U_3\left[I_{\frac{1}{3}}(\tau)\right]^{-1}$$

(55)

where $I_n(x)$ represents the modified Bessel function of the first kind. Although, Figure 4 shows significant discrepancy between Polyanin-Zaytsev solution (Equation (55)) and solution of numerical and present methods, two last ones have good sameness.

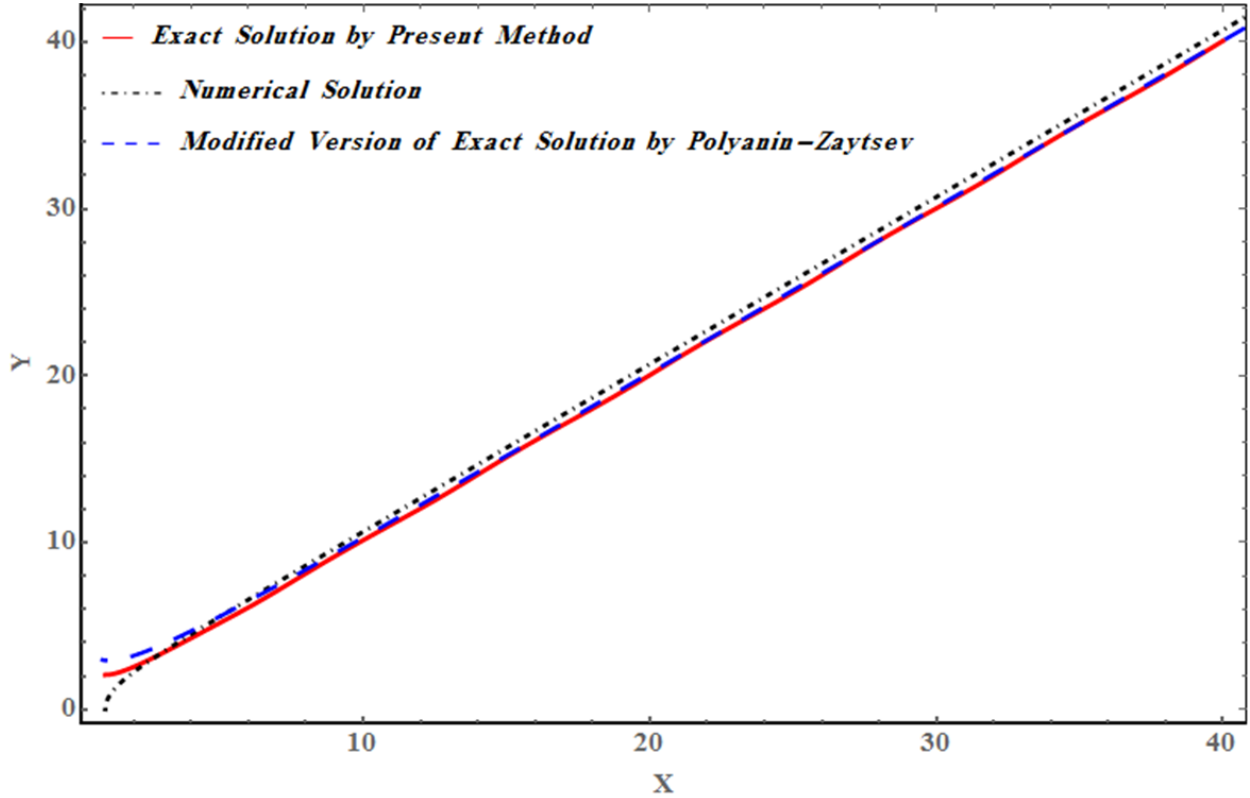

**Figure 3.** Solution of $yy'_x - y = \frac{A}{x} + \frac{A^2}{x^3}$ that A is arbitrary parameter. The Solution of present method (β=1 in Eq. (51)) has been compared by numerical solution. Also the reported exact solution of the same function that reported by [5]-Equation 1.3.1.7 requires a modification. The modified version is compared with other methods.

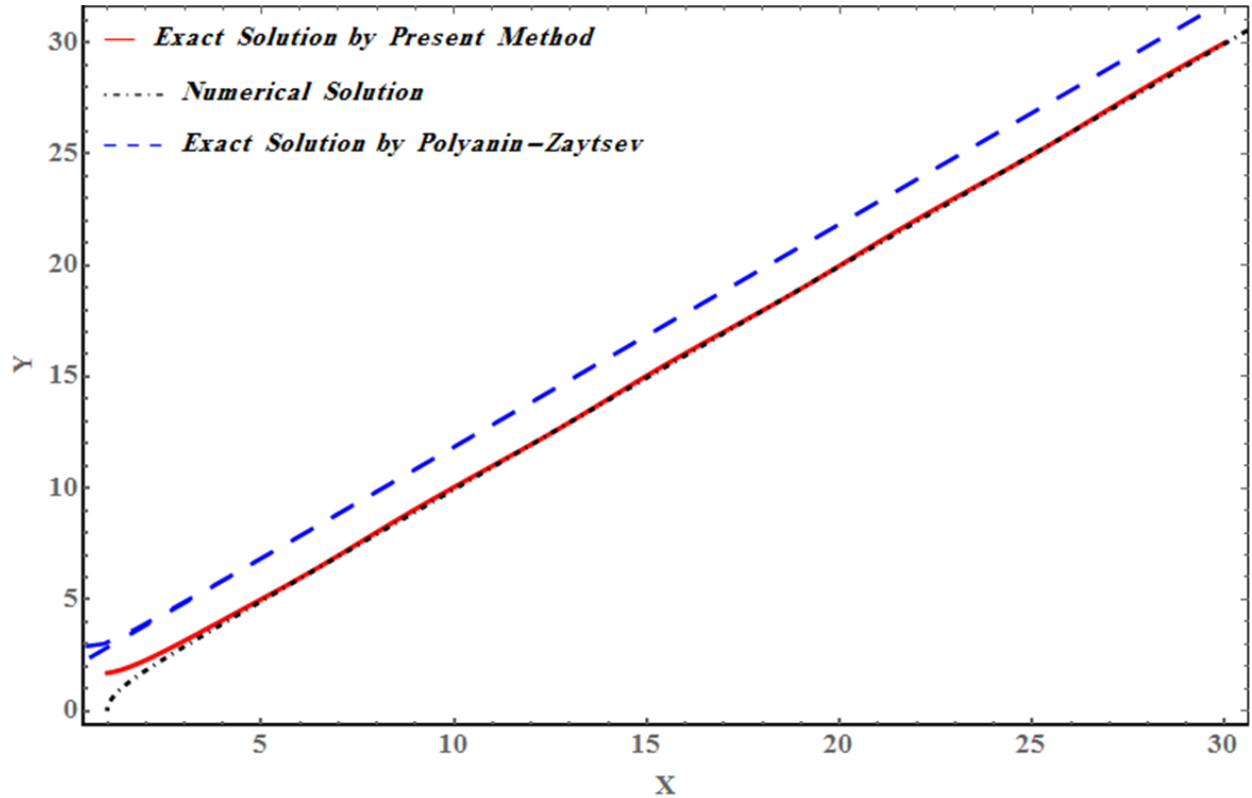

**Figure 4.** The present method (with β=1 in Eq. (51)) is compared by the numerical method and an exact solution that reported in [5](Equation 1.3.1.33) for equation $yy'_x - y = \frac{A}{x^2}$ that A is arbitrary parameter that here put 1.

### 5.3 Exponential Inhomogeneity

The third form of Polyanin-Zaytsev solution relates to nonhomogeneous part by exponential function. Hence, an exponential function is dealt with as $Q(x) = A(e^{2x/A} - 1)$ that has been given in Equation 1.3.1.9 of [5]. The solution has been presented parametrically in [5] as:

$$x = A\, Ln\left|\frac{\tau^2 + 1}{\tau}(\arctan \tau - C)\right|$$
$$y = \frac{A}{\tau}[\tau + (\tau^2 + 1)(\arctan \tau - C)]$$
(56)

The results of the present method with numerical method and also equation (56) for A=1 are brought together in Figure 5. This figure clearly approved the accuracy of present method to solve the exponential functions.

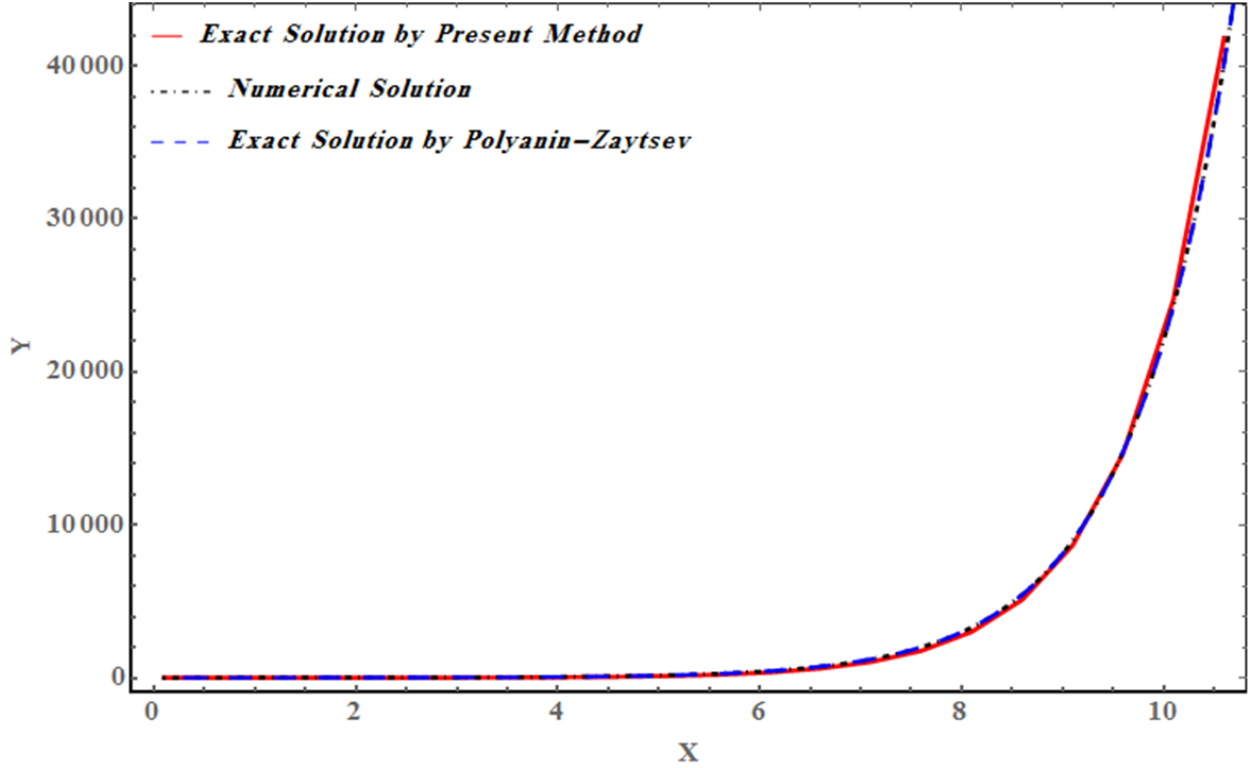

**Figure 5.** The result of the present method is compared with the numerical method and the reported exact solution by Polyanin-Zatsev ([5]-Equation 1.3.1.9). The equation is $yy'_x - y = A(e^{2x/A} - 1)$ that A is arbitrary parameter. The factor $\beta$ (Equation (51)) equals to 0.32.

## 6 Concluding Remarks

This paper has presented a reliable and exact method to solve a well-known unsolvable equation that is named Abel equation. In brief, the solution of normal form of Abel equation, $yy'_x - y = Q(x)$, is found as roots of the following cubic equation:

$$Z^3 + p\,Z + q = 0$$

$$p = -\frac{a^2}{3} + b, \quad q = 2\left(\frac{a}{3}\right)^3 - \frac{ab}{3} + c,$$

$$a = -4, \quad b = 3 - c - \frac{4Q(x)}{x + \varphi}, \tag{57}$$

$$c = \frac{0.5\,\psi\,\text{Sin}\,2x - \left[2 + \frac{1}{sgn(x)|x|}\right]\psi\,\text{Sin}^2 x + 2\psi^2\left[\text{Cos}\,x - \frac{\text{Sin}\,x}{sgn(x)|x|}\right] - \text{Sin}^3 x}{-2\,\psi^3}$$

$$\psi = sgn(x)\,|x|\,Si(x)$$

sgn(x) and Si(x) represent the sign of x and the Sine integral, respectively. Finally, the exact solution is obtained as:

$$y(x) = \frac{1}{2}(x + \varphi)\left[\beta Z(x) + \frac{1}{3}\right] \tag{58}$$

where φ is constant that is calculated by initial condition. Also, β is a factor related to the order of inhomogeneity that has been defined in Table 1.

## 8 Appendix A

Bougoffa [6] has reported a method to exact solution of the general form of the Abel equation in second kind (Equation (1)) as follow. If there exists a constant $\lambda$ such that

$$2B_1(x)\, g_0(x) = \lambda\, B_2(x)\, g_1(x) \qquad g_i(x) \neq 0,\, i = 0,1 \tag{59}$$

Then Eq. (1) admits the general solution

$$B_1(x)\, U^2 + \lambda\, B_2(x)\, U = 2 \int \frac{f_0(x)}{g_1(x)} B_1(x)\, dx + C_1 \tag{60}$$

Where:

$$B_1(x) = e^{-2 \int \frac{f_2(x)}{g_1(x)} dx} \qquad \text{and} \qquad B_2(x) = e^{-\int \frac{f_1(x)}{g_0(x)} dx} \tag{61}$$

And $C_1$ is integration constant.

## 9 Appendix B

In 1933 G. Julia introduced a functional relation between the variable coefficients of Equation (1) that can lead to the exact general solution of this equation [5, p. 27; type (b)]. Thus, if:

$$g_0[2f_2 + g'_{1_x}] = g_1[f_1 + g'_{0_x}] \qquad g_1(x) \neq 0 \tag{62}$$

Then (1) admits the general solution

$$\frac{g_1(x)\, U^2 + 2\, g_0(x)\, U}{g_1(x)\, J} = 2 \int \frac{f_0(x)}{g_1(x)\, J} dx + C_2 \quad , \quad J(x) = e^{2 \int \frac{f_2(x)}{g_1(x)} dx} \tag{63}$$

Where $C_2$ is integration constant and J is the integration factors.

## 10 Appendix C

Recall that an odd function is a function with the property that $f(-x) = -f(x)$. An odd function is said to be symmetric about the origin. When an odd function is integrated over a symmetric interval [-a, a] (a > 0) the value of the integral is zero.

$$\int_{-a}^{+a} f(x)\, dx = 0 \tag{64}$$